\documentclass[10pt,a4paper,twoside]{article}

\usepackage{a4}
\usepackage{exscale}
\usepackage[]{amsmath}
\usepackage{amssymb}
\usepackage{mathtools}
\usepackage{amsfonts}
\usepackage{eucal}
\usepackage[amsmath,thmmarks]{ntheorem}
\usepackage{bbm}

\theoremseparator{:}
\newtheorem{theorem}{Theorem}[section]

\newtheorem{proposition}[theorem]{Proposition}

\theorembodyfont{\normalfont}
\newtheorem{remark}[theorem]{Remark}
\newtheorem*{remarkon}{Remark}

\theoremsymbol{\ensuremath{\square}}
\theoremheaderfont{\normalfont\scshape}
\theorembodyfont{\normalfont}
\newtheorem*{proof}{Proof}

\newcommand{\stack}[2]{\stackrel{\substack{#1}}{#2}}

\newcommand{\R}{\ensuremath{\mathbb{R}}}
\newcommand{\N}{\ensuremath{\mathbb{N}}}

\newcommand{\ccdot}{\,\cdot\,}
\newcommand{\II}{{\mathchoice {1\mskip-4mu\mathrm l} 
     {1\mskip-4mu\mathrm l}{1\mskip-4.5mu\mathrm l} {1\mskip-5mu\mathrm l}}}

\renewcommand{\epsilon}{\varepsilon}

\usepackage[active]{srcltx}
\usepackage[hang,small,bf]{caption}
\usepackage{fancyhdr}

\usepackage[top=2.74cm,bottom=2.54cm,left=2.54cm,right=2.54cm]{geometry}
\theoremlisttype{all}
\newcommand{\comment}[1]{}
\newcommand{\longversion}[1]{}

\DeclareMathOperator{\Reg}{Reg}

\DeclareMathOperator{\diverg}{div}

\newcommand{\Ln}{\mathcal{L}^n}
\newcommand{\dx}{\, dx}

\newcommand{\kabs}[1]{\ensuremath{\vert#1\vert}}
\newcommand{\babs}[1]{\ensuremath{\big\vert#1\big\vert}}

\newcommand{\knorm}[1]{\ensuremath{\Vert#1\Vert}}

\newcommand{\tmint}{\mathop{\int\hskip -0,89em -\,}}

\usepackage{titlesec}  
\titleformat{\section}{\bf\large}{\thesection}{1em}{}
\titleformat{\subsection}{\bf\normalsize}{\thesubsection}{1em}{}

\numberwithin{equation}{section}
\allowdisplaybreaks[3]

\title{\vspace{-1cm}
{\bf Regularity versus singularity for elliptic problems in two dimensions}
}
\date{}
\author{Lisa Beck\footnote{L. Beck, Scuola Normale Superiore di Pisa, Piazza dei Cavalieri 7, Pisa, Italy. E-mail: lisa.beck@sns.it}}

\begin{document}

\vspace{-5cm}

\maketitle

\begin{abstract}
In two dimensions every solution to a nonlinear elliptic system $\diverg a(\cdot,u,Du)=0$ has H\"older continuous first derivatives provided that standard continuity, ellipticity and $p$-growth assumptions hold for some $p \geq 2$. We give an example showing that this result cannot be extended to elliptic systems in the subquadratic case, i.\,e. that weak solutions are not necessarily continuous if $1< p <2$. Furthermore, we discuss related results for variational integrals.
\end{abstract}


\section{Introduction}

The aim of this paper is to investigate some regularity properties and the possible existence of singularities for vector-valued weak solutions $u \in W^{1,p}(\Omega,\R^N)$ of second order elliptic systems in divergence form
\begin{equation}
\label{systemOmega}
- \diverg a (\, \cdot \, , u, Du ) \, = \, b(\ccdot,u,Du)  \qquad \text{in } \Omega \,.
\end{equation}
We further discuss some related results for the minimization problem of convex variational integrals
\begin{equation}
\label{minimizer}
{\cal F}[w] \, := \, \int_\Omega f(\ccdot,w,Dw) \dx
\end{equation}
in $W^{1,p}(\Omega,\R^N)$. Here the case $n,N \in \N$ for $n \geq 2$, $p \in (1,\infty)$ is considered, with $\Omega$ denoting a bounded domain in $\R^n$. The two problems are closely connected in the following sense: provided that the integrand is sufficiently regular, minimizers of ${\cal F}$ solve the Euler-Lagrange system associated to ${\cal F}$:
\begin{equation*}
\diverg D_z f(\ccdot,u,Du) \, = \, D_u f(\ccdot,u,Du) \qquad \text{ in } \Omega \,.
\end{equation*}
Nevertheless, exploiting the fact that the minimizer is a solution to the Euler-Lagrange system does often not lead to the desired results since this approach cannot distinguish between minimizers and extremals. Therefore, the regularity of weak solutions and of minimizers has to be discussed separately to a large extend. Various, by now classical results are available in the literature and helped to establish a quite general regularity theory for both the scalar ($N=1$) and the vectorial ($N>1$) case. Furthermore, several counterexamples to full regularity were constructed in the vectorial case. In what follows, we give a short description of the known regularity theory and study its  consequences, but also its limits for the two-dimensional case $n=2$. We then demonstrate how regularity and smoothness of solutions depend on the integrability exponent $p$.

We begin with a short overview on existing regularity results (for more details and an extensive list of references we recommend Mingione's invitation to the dark side \cite{Mingione06Dark}), supposing always that the coefficients or the integrands are sufficiently regular and that they satisfy suitable assumptions (see \eqref{GV-f} and \eqref{GV} below). Since the fundamental papers of De Giorgi, Nash and Moser on solutions to single equations, the theory of scalar weak solutions or minimizers is by now well understood, establishing regularity in the sense that the gradients are locally H\"older continuous, independently of the space dimension~$n$. In the vectorial case instead  first counterexamples of De\,Giorgi \cite{DEGIORGI68} and of Giusti and Miranda \cite{GIUSMIRA68EX} dating from 1968 have revealed that solutions to elliptic systems as well as minima of variational integrals may develop singularities for $n \geq 3$ even if the coefficients are analytic. Hence, in contrast to the scalar case, we can in general expect only a partial regularity result, which means regularity outside a negligible set, which is called the singular set. Here regularity is always understood as (H\"older) continuity of the solution (or of its gradient), and we introduce the set
\begin{equation*}
\Reg_{\alpha}(w) \, := \, \big\{ x \in \Omega \colon w \text{ is locally continuous with H\"older exponent } \alpha \text{ near } x \big\}
\end{equation*}
for functions $w \in L^1(\Omega,\R^N)$ and exponents $\alpha \in [0,1]$ (with the obvious inclusion $\Reg_{\alpha_1}(w) \supset \Reg_{\alpha_2}(w)$ for $\alpha_1 < \alpha_2$). Partial regularity of the solution itself in dimensions $n \leq p+2$ (which are referred to as low dimensions) is obtained via Morrey-type estimates. It is traced back to Campanato \cite{CAMPANATO82,CAMPANATO87b,CAMPANATO87} and yields in particular the bound $n-p$ on the Hausdorff dimension of the singular set. Partial regularity of the gradient instead was accomplished via Campanato-type estimates in general dimensions by various authors starting from the classical papers  \cite{MORREY68,GIUSMIRA68,GIAMOD79,EVANS86,FUSHUT86,GIAMOD86,ACEFUS87} and resulting in optimal H\"older continuity outside a set of Lebesgue measure zero. However, the counterexamples available in the literature still leave open the question of whether or not full regularity necessarily holds true in dimension $n=2$. We now discuss in more detail the main ingredients for proving regularity, namely classical Morrey- and Campanato-decay estimates for the gradient as well as some particular features exhibited in the two-dimensional case.

We are first interested in H\"older continuity of weak solutions to \eqref{systemOmega} or minimizers of \eqref{minimizer}, which will both be denoted by $u$. In view of Sobolev's embedding every function in $W^{1,q}(\Omega,\R^N)$ with $q > 2$ is continuous with some (possibly small) H\"older exponent. By taking advantage of the minimality resp. the system equation, it turns out that even if $u$ is a~priori only in $W^{1,p}(\Omega,\R^N)$ for some $p \leq 2$, then it indeed also belong to such $W^{1,q}(\Omega,\R^N)$, \emph{provided that} the a~priori integrability is not too small, i.\,e. that $p \geq p_0$ for some $p_0 \in (1,2)$ depending crucially on the structure constants. On the contrary, for small integrability exponents $p \in (1,p_0)$ only the Morrey regularity theory is available, which states the equivalence $\Reg_0(u) = \Reg_\lambda(u)$ for all $\lambda \in (0,1)$ and in fact guarantees that $\Reg_0(u)$ coincides with the whole domain $\Omega$ possibly apart from a set of Hausdorff dimension less than $2-p$ (but this does not exclude singularities/discontinuities).

In the next step a non-trivial relation between $\Reg_0(u)$ and $\Reg_0(Du)$ is established. We first recall the counterexamples \cite{NECAS77,HAOLEONEC96} of Ne\v{c}as et al., where an integrand $f$ is constructed which -- in contrast to the examples mentioned before -- depends only on the gradient variable, and where the solution to the related minimization problem \eqref{minimizer}  for dimensions $n\geq5$ (resp. its Euler-Lagrange equation for $n\geq3$) is Lipschitz-continuous, but not of class $C^1$. This example is important for two reasons: on the one hand this particular singular solution arises from the vectorial setting and not from an interaction effect with the $(x,u)$-dependency of the integrand or the coefficients; on the other hand it shows that in general the strict inclusion $\Reg_0(Du) \subset \Reg_0(u)$ holds. However, we now focus on the two-dimensional situation which is very different: in fact, $C^1$-regularity of solutions is well-known if the convex integrand resp. the coefficients of the system depend on the gradient variable, cf. Proposition~\ref{Prop-a-priori}. Moreover, by a simple comparison or perturbation argument, the regularity of the comparison solution is carried over to the solution of the original problem and implies $\Reg_0(u) = \Reg_0(Du)$.

The last step is the regularity improvement for the gradient $Du$: the minimality property of $u$ or the system equation can be used under quite general conditions (in particular for arbitrary dimension and arbitrary integrability exponents) to prove the equivalence $\Reg_0(Du) = \Reg_\beta(Du)$ with $\beta$ the optimal H\"older exponent (given in terms of the regularity of the coefficients or the integrand with respect to the $(x,u)$-variables) and to further show that this regularity criterion applies $\Ln$-almost everywhere on $\Omega$.

In conclusion, the following, straightforward strategy can be employed in two dimensions:
\begin{equation*}
\Omega \, \stack{\text{Sobolev} \\ \text{for } p \geq p_1}{=} \, \Reg_0(u) \, \stack{\text{Morrey} \\ \text{estimates}}{=} \, \Reg_\lambda(u)  \, \stack{\text{freezing} \\ \text{argument}}{=} \, \Reg_0(Du) \, \stack{\text{Campanato} \\ \text{estimates}}{=} \, \Reg_\beta(Du)
\end{equation*}
(with $\beta, \lambda \in (0,1)$ as above), and we now proceed to the precise statement of the \emph{full} regularity results for minimizers of variational integrals and for weak solutions to elliptic system: Dealing with variational integrals we consider integrands $f \colon \Omega \times \R^N \times R^{2N} \to \R$ subject to standard differentiability, growth and convexity assumptions: we require that $z \mapsto f(\cdot,\cdot,z)$ is of class $C^2$ with jointly continuous second order derivatives and that we have
\begin{equation}
\label{GV-f}
\left\{ \quad \begin{array}{l}
D_{zz} f(x,u,z) \text{ is continuous on } \Omega \times \R^N \times R^{2N} \,, \\[0.1cm]
\nu \, \kabs{z}^p \, \leq \, f(x,u,z) \, \leq \, L \, (1+\kabs{z})^p \,, \\[0.1cm]
\nu \,  (1+\kabs{z})^{p-2} \kabs{\lambda}^2  \, \leq \, 
D_{zz} f(x,u,z) \, (\lambda , \lambda) \, \leq \, L \,  (1+\kabs{z})^{p-2} \kabs{\lambda}^2 \,, \\[0.1cm]
\babs{D_zf(x,u,z) -  D_zf(\bar{x}, \bar{u}, z)} \, \leq \, L \, (1 + 
	\kabs{z})^{p-1} \, \omega_{\alpha_1}\big(\kabs{x-\bar{x}} + \kabs{u-\bar{u}}\big)\,, \\[0.1cm]
\babs{f(x,u,z) -  f(x, \bar{u}, z)} \, \leq \, L \, (1 + 
	\kabs{z})^{p} \, \omega_{\alpha_2}\big(\kabs{u-\bar{u}}\big) \,,
\end{array} \right.
\end{equation}
for all $x,\bar{x} \in \Omega$, $u, \bar{u} \in \R^N$, $z,\lambda \in \R^{2N}$, and with fixed $L \geq \nu > 0$, $\alpha_1, \alpha_2 \in (0,1)$. Here $\omega_{\beta} \colon \R^+ \to \R^+$ denotes for arbitrary $\beta \in (0,1]$ the modulus of continuity
\begin{equation*}
\omega_{\beta}(t) \, = \, \min \{1,t^{\beta}\}. 
\end{equation*}
Then the following full regularity result holds:

\begin{theorem}
\label{reg_varint}
Let $p \in (1,\infty)$ and suppose that $\Omega \subset \R^2$ is a bounded domain. There exists $p_0 = p_0(N,\nu,L)<2$ such that the following statement is true: whenever $u \in W^{1,p}(\Omega,\R^N)$ is a minimizer to \eqref{minimizer} under the assumptions \eqref{GV-f}, then $Du$ is locally H\"older continuous in the interior of $\Omega$ with optimal exponent $\beta := \min \{ \frac{\alpha_2}{2-\alpha_2}, \frac{\alpha_1}{2} \} < \frac{1}{2}$, i.\,e. $u \in C^{1,\beta}_{\rm{loc}}(\Omega,\R^N)$, for $p \geq p_0$. The same assertion remains true for all $p > 1$ if the integrand is independent of $u$, i.\,e. $f(x,u,z) \equiv f(x,z)$, or if $u \in C^{0}_{\rm{loc}}(\Omega,\R^N)$.
\end{theorem}

For the treatment of elliptic systems we consider coefficients $a \colon \Omega \times \R^N \times \R^{2N} \to \R^{2N}$ for which we impose similar assumptions concerning differentiability, growth and ellipticity: we require that we have
\begin{equation}
\label{GV}
\left\{ \quad \begin{array}{l}
z \mapsto a(x,u,z) \text{ is of class } C^1(\R^{2N},\R^{2N})  \,, \\[0.1cm]
\babs{a(x,u,z)} + \babs{ D_z a(x,u,z) } \, \big(1 + \kabs{z}\big) \, \leq \, L \, \big(1 + \kabs{z} \big)^{p-1} \,, \\[0.1cm]
D_z a (x,u,z) \, \lambda \cdot \lambda \, \geq
    \, \nu \, \big(1 + \kabs{z} \big)^{p-2} \kabs{\lambda}^2 \,, \\[0.1cm]
 \babs{a(x,u,z) -  a(\bar{x}, \bar{u}, z)} \, \leq \, L \, \big(1 +
    \kabs{z} \big)^{p-1} \, 
	\omega_{\alpha}\big(\kabs{x-\bar{x}} + \kabs{u-\bar{u}}\big) \,,
\end{array} \right.
\end{equation}
for all $x,\bar{x} \in \Omega$, $u, \bar{u} \in \R^N$, and $z,\lambda \in \R^{2N}$, with $\alpha \in (0,1)$. For the inhomogeneity $b \colon \Omega \times \R^N \times \R^{2N} \to \R^{N}$ we assume the controllable growth condition
\begin{equation}
\label{inhomo-growth}
\kabs{b(x,u,z)} \, \leq \, L (1 + \kabs{z})^{p-1}
\end{equation} 
for all $x \in \Omega$, $u \in \R^N$, and $z \in \R^{2N}$. The corresponding regularity result is then given as follows: 

\begin{theorem}
\label{reg_systems}
Let $p \in (1,\infty)$ and suppose that $\Omega \subset \R^2$ is a bounded domain. There exists $p_1 = p_1(N,\nu,L) < 2$ such that the following statement is true: whenever $u \in W^{1,p}(\Omega,\R^N)$ with $p \geq p_1$  is a weak solution to \eqref{systemOmega} under the assumptions \eqref{GV} and \eqref{inhomo-growth}, then $Du$ is locally H\"older continuous in the interior of $\Omega$ with optimal exponent $\alpha$, i.\,e. $u \in C^{1,\alpha}_{\rm{loc}}(\Omega,\R^N)$. The same assertion remains true for all $p > 1$ if the coefficients are independent of $u$, i.\,e. $a(x,u,z) \equiv a(x,z)$, or if $u \in C^{0}_{\rm{loc}}(\Omega,\R^N)$.
\end{theorem}

\begin{remarkon}
This regularity statement is extended easily to bounded solutions of inhomogeneous systems under a natural growth condition and the additional standard smallness assumption on $\knorm{u}_{L^{\infty}}$.
\end{remarkon}

The first regularity result Theorem~\ref{reg_varint} is a special case of \cite[Theorem~1.4]{BECK09} (using in turn the arguments from \cite[Theorem~1.7]{KRIMIN05}) and is recalled here in order to draw a picture, as complete as possible, of the topic of regularity for two-dimensional elliptic problems. The second result seems not to be stated explicitly in the literature -- even if all parts of the proof of Theorem~\ref{reg_systems} are essentially known. For this reason we give a proof in Section~\ref{sec-reg}. For the optimality in both theorem we refer to \cite[Example~1.1]{GROTOWSKI02} and \cite[Section~I]{PHILLIPS83}.

We now comment on the existing literature and the reason for which the distinction concerning the $u$-dependence in the above statements is made. For the moment we shall restrict ourselves to the $u$-independent case and concentrate on systems rather than on variational integrals (for which we can pass to the Euler-Lagrange system). Imposing a differentiable dependence on the $x$-variable, Star\'a \cite{STARA71} succeeded in showing the existence of higher order derivatives, ending up with a global H\"older continuity result for $Du$. In the case of merely H\"older continuous dependence on $x$, arguments of Campanato \cite[Section~3]{CAMPANATO82} reveal that every solution $u \in W^{1,p}(\Omega,\R^N)$ to \eqref{systemOmega} is in fact H\"older continuous independently of the value of $p \in (1,\infty)$. This corresponds to the first step described in the strategy above. Similarly, regarding fractional Sobolev spaces as a generalization of the class of H\"older continuous functions, one further has a fractional differentiability result, see \cite{MINGIONE03arch}. 

For the $u$-dependent case and with an a~priori continuous solution $u$, the low dimensional theory guarantees in a first step H\"older continuity of $u$. At this point the philosophy is to pass to new coefficients $\tilde{a}(x,z) = a(x,u(x),z)$, and one is then back in the $u$-independent case.

After having given the background for the full regularity results, we now proceed to the main objective of this paper: we want to address the problem of whether full regularity necessarily holds for all $p>1$ (which would be equivalent to setting $p_0,p_1 = 1$ in Theorem~\ref{reg_varint}, Theorem~\ref{reg_systems} above) or whether there might arise singularities -- a question which was posed by Campanato \cite{CAMPANATO82}.


With the strategy from above, the optimal H\"older regularity of $Du$ follows only outside of the open set $\Reg_0(u)$, i.\,e. in the case $p \in (1,p_1)$ outside a negligible set of Hausdorff dimension less than $2-p$, and it is not clear to what extend this result can still be improved. In light of the full regularity results, the construction of an adequate system resp. functional that might provide an example with a singular solution demands one of the following features: either the integrand resp. coefficients have to be less regular in the $x$-variable (not satisfying \eqref{GV-f}$_4$ resp. \eqref{GV}$_4$), or they have to depend explicitly on the $u$-variable (in case of continuous dependence we also need a construction involving a discontinuous solution). We shall deal with both situations, working with the function $u(x)= x/\kabs{x} \in W^{1,p}(B^n,\R^n)$ for all $p \in (1,n)$, which is discontinuous in one point and appears as a prominent example in the literature: Giusti and Miranda \cite{GIUSMIRA68EX} constructed a quadratic-type functional which is minimized by $u$ for $n \geq 3$. Moreover, passing to the related Euler-Lagrange equation, it is at the same time also a weak solution to an elliptic system. Taking advantage of this construction we will provide in Section~\ref{sec_irr} for $n=2$ and every $p \in (1,2)$ an example of a functional and a system with $L^{\infty}$-dependence on $x$ which are solved by $u$. In the next step we will rewrite the coefficients, ending up with:

\begin{theorem}
\label{existence-example}
Let $u \colon \R^2 \supset B \to \R^2$ be given by $u(x) = x/\kabs{x}$ and let $p \in (1,2)$. Then $u \in W^{1,p}(B,\R^2) \cap L^{\infty}(B,\R^2)$, and there exist coefficients $a \colon \R^N \times \R^{2N} \to \R^{2N}$ satisfying the assumptions \eqref{GV} for some $L \geq \nu > 0$ and every $\alpha \in (0,1)$ such that $u$ is a weak solution of the homogeneous system $\diverg a(u,Du) = 0$ in $B$.
\end{theorem}

Theorem~\ref{existence-example} is established via Proposition~\ref{prop-ass-coefficients} and  Theorem~\ref{theorem-ex}, and it further yields the strict inequality $p_1 > 1$. At this point we recall that $p_1$ crucially depends on the structure constants, in particular on the ratio $L/\nu$. As a consequence the closer the integrability exponent $p$ is chosen to $2$, the greater this ratio needs to be chosen; see Remark~\ref{rem-blow-up}. However, whether or not there exists a minimization problem with singular solution under the assumptions \eqref{GV-f} remains open.

We mention briefly that $x/\kabs{x}$ is frequently considered as a function taking almost everywhere values in the unit-sphere $S^{n-1} \subset \R^n$, i.\,e. as a function in the space $W^{1,p}(B^n,S^{n-1})$ for $p<n$. Direct computation shows that $x/\kabs{x}$ is $p$-harmonic in the sense that it satisfies
\begin{equation*}
\int_{B^n} \kabs{Du}^{p-2} Du \cdot D\varphi \dx \, = \, \int_{B^n} \kabs{Du}^p u \cdot \varphi \dx
\end{equation*}
for all $\varphi \in W_0^{1,p}(B^n,\R^n) \cap L^{\infty}(B^n,\R^n)$. In other words, it is stationary for the $p$-energy $\int_{B^n} \kabs{Dw}^p \dx$ under the constraint $\kabs{w}=1$ almost everywhere (note that this case is not contained in the above regularity theory, but it is covered by a serial of classical papers). In fact, even minimality holds, see \cite{HARLINWAN98,HONG01}.

We close the introduction with some remarks on the notation: we write $B_{\rho}(x_0) := \{x \in \R^n: \kabs{x-x_0} < \rho\} $ for the $n$-dimensional ball centered at $x_0 \in \R^n$ with radius $\rho > 0$. The function spaces used in this paper are the H\"older spaces $C^{k,\alpha}$, Morrey-spaces $L^{p,\sigma}$, the Sobolev spaces $L^{p}$ and (fractional) Sobolev spaces $W^{\theta,p}$, with $\alpha,\theta \in (0,1]$, $k \in \N$, $\sigma > 0$ and $p \in [1,\infty)$ (see e.g. \cite[Chapter 7]{ADAMS75} for the definition and embedding theorems for fractional Sobolev spaces). Furthermore, we shall use two abbreviations: for a bounded set $X \in \R^n$ with positive Lebesgue-measure we denote the average of a function $f \in L^1(X)$ by $\tmint_{X} f \dx$, and for $\xi \in \R^k$ we write $V(\xi) = (1+\kabs{\xi}^2)^{(p-2)/4} \xi$.


\section{Review of some regularity results}
\label{sec-reg}

We collect some regularity results which are partially already available in the literature. We will only outline the proofs or give suitable references. We first observe that every 
weak solution $u \in W^{1,p}(\Omega,\R^N)$ to \eqref{systemOmega} with coefficients not depending explicitly on $u$ (or with $u$ being a priori H\"older continuous) actually belongs to a fractional Sobolev space. More precisely, we have 

\begin{proposition}
\label{Lemma-frac}
Let $u \in W^{1,p}(\Omega,\R^N)$, $p \in (1,\infty)$, be a weak solution to \eqref{systemOmega} under the assumption \eqref{GV} and \eqref{inhomo-growth}. Furthermore, we suppose the coefficients to be independent of $u$, i.\,e. $a(x,u,z) \equiv a(x,z)$, or $u \in C^{0,\gamma}_{\rm{loc}}(\Omega,\R^N)$ for some $\gamma > 0$. Then $V(Du) \in W^{\alpha',2}_{\rm{loc}}(\Omega,\R^N)$ and $Du \in W^{\min\{1,2/p\}\alpha',p}_{\rm{loc}}(\Omega,\R^N)$ for every $\alpha' < \alpha$.

\begin{proof}[Sketch]
Using difference quotients techniques we can argue similarly to Mingione \cite[proof of Proposition 3.1]{MINGIONE03arch} and \cite[proof of Proposition 5.2]{MINGIONE03}, where the superquadratic case was considered, in order to derive the fractional differentiability from an (uniform) estimate for finite difference quotients. It should be noted that no further assumption (such as the continuity assumption \cite[(1.8)]{MINGIONE03arch}) on the inhomogeneity is needed.
\end{proof}
\end{proposition}

After this higher differentiability result (which implies higher integrability via fractional Sobolev embedding), we come to an essential ingredient needed for the application of the comparison argument in the proof of Theorem~\ref{reg_systems}, namely a~priori estimates for solutions of a ``frozen'' problem. Following \cite[Section 3]{CAMPANATO82} we see that these solutions admit second order derivatives. Using Gehring's lemma in order to deduce a higher integrability result of second order derivatives (or applying a version of Widman's hole filling trick \cite{WIDMAN71} as in \cite[Lemma 8.2]{SCHEVENSCHMIDT10}), we thus obtain (see also \cite[Theorem~3.I]{CAMPANATO87} for the superquadratic case):

\begin{proposition}
\label{Prop-a-priori}
Let $v \in W^{1,p}(B_R(x_0),\R^N)$ be a weak solution to 
\begin{equation*}
\diverg a_0(Dv) = 0  \qquad \text{in } B_R(x_0) \subset \Omega \subset \R^n
\end{equation*}
with coefficients $a_0(\cdot)$ under the assumptions \eqref{GV}$_{0,1,2}$. Then there exists $0 < \epsilon = \epsilon(n,N,p,L,q)$ such that for every $\rho \in (0,R]$ we have:
\begin{align*}
& \int_{B_{\rho}(x_0)} \babs{V(Dv)-\big(V(Dv)\big)_{B_{\rho}(x_0)}}^2 \dx \, \leq \, c \, \Big( \frac{\rho}{R}\Big)^{2 + \epsilon}
	\int_{B_{R}(x_0)} \babs{V(Dv)}^2 \dx
\shortintertext{and}
& \int_{B_{\rho}(x_0)} \babs{V(Dv)}^2 \dx \, \leq \, c \, \Big( \frac{\rho}{R}\Big)^{\min\{n,2+\epsilon\}} 
	\int_{B_{R}(x_0)} \babs{V(Dv)}^2 \dx \,,
\end{align*}
and both constants depend only on $n,N,p,L$ and $\nu$.
\end{proposition}

This proposition uncovers a peculiarity of the two-dimensional case $n=2$: the solution to the comparison problem has H\"older continuous first derivatives. This helps us to obtain in a first step a Morrey-space regularity result for the weak solution $u$ to \eqref{systemOmega} (see also \cite[proof of Theorem 1.I]{CAMPANATO82}), which in turn yields the global H\"older regularity result in the interior of $\Omega \subset \R^2$:

\begin{proof}[of Theorem~\ref{reg_systems}]
We here follow the arguments in \cite[Section 9]{KRIMIN05} where the related result for minimizers of variational functionals was proved for $p \geq 2$.

\emph{Step 1a: Determination of $p_1$ and preliminary regularity improvement of $u$.} Via a standard Caccioppoli inequality and Gehring's Lemma we first recall the higher integrability result $Du \in L^q_{\rm{loc}}(\Omega,\R^N)$ for some exponent $q>p$ depending only on $N,p$ and $\frac{L}{\nu}$. This is exploited to determine the number $p_1(N,L,\nu)<2$ such that $q(N,p,L,\nu)>2$ for all $p \geq p_1$. Sobolev's embedding in turn implies $u \in C^{0,\lambda}_{\rm{loc}}(\Omega,\R^N)$ for some $\lambda = \lambda(N,L,\nu)>0$. For $p \in (1,p_1)$ and general coefficients instead, local continuity of $u$ is assumed. In fact, this is equivalent to local H\"older continuity $u \in C^{0,\lambda}_{\rm{loc}}(\Omega,\R^N)$ for some $\lambda>0$ by the low dimensional theory, see \cite[Theorem 1.I]{CAMPANATO82}.

\emph{Step 1b: Morrey-space regularity: $Du \in L^{p,2-\mu}_{\rm{loc}}(\Omega,\R^{2N})$ for every $\mu > 0$.} Let $B_{2R}(x_0) \subset \Omega$ and define $v \in u + W^{1,p}_0(B_R(x_0),\R^N)$ as the unique solution to
$\diverg a_0(Dv) = 0$ in $B_R(x_0)$,
where the coefficients are defined by freezing via $a_0(z) := a(x_0,(u)_{B_R(x_0)},z)$ (we note that existence and uniqueness follow from standard theory for monotone operators). Using $u-v \in W^{1,p}_0(B_R(x_0),\R^N)$ as a test-function in the weak formulation of the comparison Dirichlet problem, we deduce the energy estimate
\begin{equation*}
\int_{B_R(x_0)} \kabs{V(Dv)}^2 \dx \, \leq \, c(p,L,\nu) \, \int_{B_R(x_0)} \big( 1 + \kabs{V(Du)}^2 \big) \dx \,.
\end{equation*}
Furthermore, taking into account the growth, ellipticity and continuity assumption in \eqref{GV} as well as Poincar\'{e}'s inequality, we find the comparison estimate
\begin{align}
\label{comp-est}
\lefteqn{ \hspace{-0.5cm} c^{-1}(N,p,\nu) \int_{B_R(x_0)} \babs{V(Du) - V(Dv)}^2 \dx 
	\, \leq \, \int_{B_R(x_0)} \big[ a_0(Du) - a_0(Dv) \big] \, (Du - Dv) \dx } \nonumber \\
	& = \, \int_{B_R(x_0)} \big[ a_0(Du) - a(x,u,Du) \big] \, (Du - Dv) \dx 
	+ \int_{B_R(x_0)} b(x,u,Du) \, (u-v) \dx\nonumber \\
	& \leq \, c(N,p,L) \, R^{\alpha \lambda} \int_{B_R(x_0)} \big( 1 + \kabs{V(Du)}^2 \big) \dx \,.
\end{align}
To obtain the last line, different cases need to be distinguished: on the one hand we might be concerned with coefficients with explicit $u$-dependencies and the prerequisite $p > p_1$ or $u$ a~priori continuous. This situation is handled via the local $C^{0,\lambda}$-regularity of $u$ according to Step 1a. On the other hand the coefficients might have no explicit $u$-dependency, for which the above estimate holds trivially with $\lambda = 1$.
Combining the previous two estimates with the decay estimate in Proposition~\ref{Prop-a-priori}, we hence end up with
\begin{align*}
\int_{B_{\rho}(x_0)} \big(  1 + \kabs{V(Du)}^2 \big) \dx 
	& \leq \, c(N,p,L,\nu) \, \Big( \Big( \frac{\rho}{R}\Big)^2 + R^{\alpha \lambda} \Big) 
	\int_{B_R(x_0)} \big( 1 + \kabs{V(Du)}^2 \big) \dx 
\end{align*}
for all $\rho \in (0,R]$. An iteration procedure (see e.g. \cite[Chapter III, Lemma 2.1]{GIAQUINTA83}) then yields: for every $\mu>0$ there exists a radius $R_0 = R_0(\mu,N,p,L,\nu,\alpha) >0$ (independent of the center $x_0$ of the balls) such that
\begin{equation*}
\int_{B_{\rho}(x_0)} \big(  1 + \kabs{V(Du)}^2 \big) \dx \, \leq \, \rho^{2-\mu}
\end{equation*}
for all $\rho < R_0$. To conclude the desired Morrey-space embedding for $Du$ we observe that for large radii $\rho \geq R_0$ the left-hand side is easily estimated by $c(\mu,n,N,p,L,\nu,\alpha,\knorm{Du}_{L^p(B,\R^N)}) \rho^{2-\mu}$. We lastly note that this Morrey-type estimate is in particular a further regularity improvement of $u$ in the sense that $u$ is locally H\"older continuous for any exponent less than $1$.

\emph{Step 2: Continuity of $Du$.} We next apply the comparison estimate \eqref{comp-est} and the a priori estimate for $v$ from Proposition~\ref{Prop-a-priori}, and we find
\begin{align*}
\lefteqn{ \hspace{-0.5cm} \int_{B_{\rho}(x_0)} \babs{V(Du)-\big(V(Du)\big)_{B_{\rho}(x_0)}}^2 \dx } \\
	& \leq \, c \, R^{\alpha \lambda} \int_{B_R(x_0)} \big( 1 + \kabs{V(Du)}^2 \big) \dx
	+ c \, \Big( \frac{\rho}{R}\Big)^{2 + \epsilon} 
	\int_{B_{R}(x_0)} \babs{V(Dv)}^2 \dx \\
	& \leq \, c\big(\mu,N,p,L,\nu,\alpha,\knorm{Du}_{L^p(B,\R^N)}\big) \, 
	\big[ R^{2+ \alpha \lambda + \epsilon} + \rho^{2+\epsilon} \big] \, R^{-\mu-\epsilon }
\end{align*}
for every $\mu > 0$. Now we choose $R$ as a power of $\rho$ such that all powers of $\rho$ on the right-hand side are equal, i.\,e. $R = \rho^{(2 + \epsilon)/(2 + \alpha \lambda + \epsilon )} > \rho$. Hence, we get
\begin{equation*}
\int_{B_{\rho}(x_0)} \babs{V(Du)-\big(V(Du)\big)_{B_{\rho}(x_0)}}^2 \dx 
	\, \leq \, c \, 
	\rho^{(2+\epsilon) \frac{2 + \alpha \lambda - \mu}{2 + \alpha \lambda + \epsilon }}
\end{equation*}
for a constant $c$ admitting the same dependencies as above. Therefore, the exponent at the right-hand side is strictly greater than the space dimension $n=2$ if we choose $\mu \in(0, \frac{\epsilon \alpha \lambda}{\epsilon + 2})$ sufficiently small. Since the estimate is independent of the ball under consideration, we conclude from Campanato's characterization of H\"older continuous functions that $Du$ is in particular locally continuous in the interior of $\Omega$. 

\emph{Step 3: Optimal H\"older regularity of $Du$.} Standard regularity regularity (e.\,g. summarized for all possible exponents in \cite{Mingione06Dark} in Theorem~4.4 and the following characterization of the singular set) may now be applied, which states that local continuity of $Du$ is in fact equivalent to local H\"older continuity of $Du$ with optimal H\"older exponent $\alpha$. We thus get the desired regularity $Du \in C^{0,\alpha}_{\rm{loc}}(\Omega,\R^{2N})$.
\end{proof}

\begin{remarkon}
In fact, also global regularity estimates can be achieved in a similar way. For this purpose one supposes that $\Omega$ is a domain of class $C^{1,\alpha}$ and then studies solutions in the space $g + W^{1,p}_0(\Omega,\R^N)$ with $g \in C^{1,\alpha}(\bar{\Omega},\R^N)$. Via a standard flattening and transformation procedure the problem is first reduced to the model situation of the unit half-ball and zero-boundary values. Then all the results above need to be extended up to the boundary: for the extension of the fractional Sobolev estimates in Proposition~\ref{Lemma-frac} we refer to the approach in \cite[Proposition 5.1]{DUZKRIMIN05}, for the a~priori estimates of the frozen solution in order to conclude the H\"older regularity of $Du$ to \cite[Section 6]{CAMPANATO87} and \cite[Section 3]{BECK09b}. A (quite technical) combination of the interior and the boundary estimates then yields the global result.
\end{remarkon}


\section{An example for irregularity}
\label{sec_irr}

As already explained in the introduction, the previous regularity results still leave open the questions, namely whether there exist systems and variational integrals in the subquadratic case $p \in (1,2)$ which admit discontinuous weak solutions. The construction of such integrands or system -- depending discontinuously on the $x$-variable or depending explicitly on the solution $u$ -- shall be addressed to in this last section. Giusti and Miranda \cite{GIUSMIRA68EX} succeeded in showing that in dimensions $n \geq 3$ the function $x/\kabs{x} \in W^{1,2}(B^3,\R^3) \cap L^{\infty}(B^3,\R^3)$ is a minimizer to a quadratic-type functional and weak solution to a quasi-linear elliptic system, and they thus demonstrated that discontinuities may occur in dimensions $n\geq 3$. We now take advantage of their construction and show that in the two-dimensional case for any $p \in (1,2)$ the map $x/\kabs{x} \in W^{1,p}(B^2,\R^2)$ is a minimizer of functional satisfying the $p$-growth conditions \eqref{GV-f}, but discontinuous in the $x$-variable, and it is also a weak solution of a homogeneous elliptic system satisfying all assumptions in \eqref{GV}.

Analogously to Giusti's and Miranda's construction \cite{GIUSMIRA68EX} we start by defining a bilinear form on $\R^{2 \times 2}$ via
\begin{equation*}
A^{\kappa \lambda}_{ij}(u) \, = \, \delta_{\kappa \lambda} \delta_{ij} + \Big( \delta_{\kappa i} + \frac{2p}{2-p} \, \frac{u_i u_{\kappa}}{1+\kabs{u}^2} \Big) \, \Big( \delta_{\lambda j} + \frac{2p}{2-p} \, \frac{u_j u_{\lambda}}{1+\kabs{u}^2} \Big)
\end{equation*} 
for $p \in (1,2)$, all $u \in \R^2$ and indices $\kappa,\lambda,i,j \in \{1,2\}$. In what follows we shall use the convention $A(u)(z,\bar{z}) = \sum_{\kappa,\lambda,i,j \in \{1,2\}} A^{\kappa \lambda}_{ij}(u) z_i^{\kappa} \bar{z}_j^{\lambda}$ for all $z,\bar{z} \in R^{2 \times 2}$. We further introduce the abbreviations
\begin{align*}
T_u(z) \, := \, {\rm Tr} (z) + \frac{2p}{2-p} \, \frac{z \cdot u \otimes u}{1+\kabs{u}^2} \,.
\end{align*}
Due to the symmetry of $A$ we immediately find the following two useful identities
\begin{align}
A(u) z & = \, z + T_u(z) \, \Big( \II + \frac{2p}{2-p} \, 
	\frac{u \otimes u}{1 + \kabs{u}^2} \Big) \,, \nonumber \\
\label{Prop-ident-2}
A(u) (z ,\bar{z}) & = \, z \cdot \bar{z} +  T_u(z) \, T_u(\bar{z}) \,.
\end{align}

We next take $g:\R \to [0,1]$ as a symmetric, smooth cut-off function satisfying $\II_{\{0\}} \leq g \leq \II_{(-1,1)}$. We set 
\begin{equation*}
m_g \, := \, (p-1)^{-1}\, \big( 1 + \sup_{s \in \R} \{\kabs{g'(s)} + 2 \kabs{g''(s)} s \} \big) \, > \, 1
\end{equation*}
(the only benefit of this constant will be to compensate the effects of $g$ occurring in the convexity condition). We then define an integrand $f(x,z) \colon \R^2 \times \R^{2 \times 2} \to \R$ via
\begin{equation}
\label{integrand}
f(x,z) \, := \, \big( g(\kabs{z}^2) +  m_g \, A(x/\kabs{x})(z,z) \big)^{\frac{p}{2}}
\end{equation}
for all $x \in \R^2 \setminus \{0\} \times \R^{2 \times 2}$ (and with arbitrary value for $x=0$). By definition, the integrand is bounded (for fixed $z$) and  $0$-homogeneous in the $x$-variable, and it also satisfies the subquadratic growth and ellipticity condition:

\begin{proposition}
\label{prop-GV-f}
The integrand $f(\cdot,\cdot)$ defined in \eqref{integrand} is smooth with respect to the variable $z$ and satisfies the assumptions \eqref{GV-f}$_1$ -- \eqref{GV-f}$_3$ with constants $\nu,L$ depending only on $p,m_g$.

\begin{proof}
The smoothness of the integrand with respect to the gradient variable is guaranteed by construction. The assumption \eqref{GV-f}$_2$ on coercivity and boundedness is easily verified by taking into account the identity \eqref{Prop-ident-2}, for some constant $\nu$ and $L$ depending only on $p$ and $m_g$ (via the choice of $g$). Here we already note that $L$ blows up as $p \nearrow 2$. Hence, it only remains to prove \eqref{GV-f}$_3$: for this purpose we first observe
\begin{align*}
D_{zz} f(x,z) (\lambda,\lambda) 
	& = \, p \, \big( g(\kabs{z}^2) + m_g
	\, A(x/\kabs{x}) (z,z) \big)^{\frac{p-4}{2}}  \\
	& \quad \Big[ \big(g(\kabs{z}^2) + m_g
	\, A(x/\kabs{x}) (z,z) \big) 
	\, \big( g'(\kabs{z}^2) \, \kabs{\lambda}^2 + 2 
	\, g''(\kabs{z}^2) \, (z \cdot \lambda)^2 + m_g \, A(x/\kabs{x}) (\lambda,\lambda) \big) \\
	& \qquad {}- (2-p) \, \big( g'(\kabs{z}^2) \, z \cdot \lambda 
	+ m_g \, A(x/\kabs{x}) (z,\lambda)\big)^2 \Big]
\end{align*}
for all $x \in \R^2 \setminus \{0\}$ and $z,\lambda \in R^{2 \times 2}$. The second inequality in \eqref{GV-f}$_3$ then follows immediately from \eqref{Prop-ident-2}, whereas for the first one we need to take advantage of the Cauchy-Schwarz inequality, of $g(s) \leq 0$ for all $s \geq 0$, and of the definition of the constant $m_g$ to infer:
\begin{align*}
D_{zz} f(x,z) (\lambda,\lambda) 
	& \geq \,  p \, \big( g(\kabs{z}^2) + m_g
	\, A(x/\kabs{x}) (z,z) \big)^{\frac{p-4}{2}}  \\
	& \quad \Big[ \big(g(\kabs{z}^2) + m_g
	\, A(x/\kabs{x}) (z,z) \big) 
	\, \big( (m_g - g'(\kabs{z}^2) - 2 \, g''(\kabs{z}^2) \, \kabs{z}^2) 
	\, A(x/\kabs{x}) (\lambda,\lambda) \big) \\
	& \qquad {}- (2-p) \, m_g^2 \, A(x/\kabs{x}) (z,z) \, A(x/\kabs{x}) (\lambda,\lambda) \Big] \\
	& \geq \,  p \, \big( g(\kabs{z}^2) + m_g
	\, A(x/\kabs{x}) (z,z) \big)^{\frac{p-2}{2}} \, A(x/\kabs{x}) (\lambda,\lambda) 
	\, \geq \, c^{-1}(p,m_g) \, (1 + \kabs{z}^2)^{\frac{p-2}{2}} \, \kabs{\lambda}^2 \,.
\end{align*}
This completes the proof of the proposition.
\end{proof}
\end{proposition}

The regularity of $f$ also allows us to study the Euler-Lagrange system for ${\cal F}[\cdot]$ with integrand \eqref{integrand}, which is given by
\begin{equation}
\label{EL-Ex}
\diverg \big[ \big( g(\kabs{Du}^2) + m_g \, A(x/\kabs{x})(Du,Du) \big)^{\frac{p-2}{2}} \, \big( g'(\kabs{z}^2) \, Du + m_g \, A(x/\kabs{x}) Du \big) \big] \, = \, 0 \,.
\end{equation}
Due to the convexity of $f$ minimizers and critical points of the functional \eqref{minimizer} indeed coincide. This fact is now exploited to determine a discontinuous minimizer $u$, which in particular demonstrates that both the singular sets of $u$ and of $Du$ are not empty. Moreover, by the strict convexity $u$ is even the unique minimizer with respect to its own boundary values.

\begin{proposition}
\label{prop-min-Ex}
Assume $\Omega = B \subset \R^2$, $p \in (1,2)$, and let $u \colon B \to \R^2$ be given by $u(x) = x/\kabs{x}$. Then $u \in W^{1,p}(B,\R^2) \cap L^{\infty}(B,\R^2)$, and $u$ is the unique minimizer of the functional \eqref{minimizer} with integrand $f(\cdot,\cdot)$ defined in \eqref{integrand} among all functions in the class $u + W_0^{1,p}(B,\R^2)$.

\begin{proof}
We start with some preliminary observations and calculations: we note that $u$ is smooth in $\R^2 \setminus \{0\}$ and $\kabs{u}$ is bounded by $1$. Furthermore, for all $x \in \R^2 \setminus \{0\}$ and every $\kappa \in \{1,2\}$ we find
\begin{equation*}
2 \, \sum_{i \in \{1,2\}} u_i D_{\kappa} u_i \, = \, D_{\kappa} \kabs{u}^2 \, = \, 0 \,.
\end{equation*}
We next calculate for all $x \in \R^2 \setminus \{0\}$
\begin{align*}
Du(x) & = \, \frac{\II}{\kabs{x}} - \frac{x \otimes x}{\kabs{x}^3} \qquad \text{with } \kabs{Du} \, = \, {\rm Tr} (Du) \, = \, \kabs{x}^{-1} \,, \\
A(x/\kabs{x}) Du & = \, \frac{2 \, \II}{\kabs{x}} + \frac{2(p-1)}{(2-p)} \, \frac{x \otimes x}{\kabs{x}^3} \,, \\[0.1cm]
A(x/\kabs{x}) (Du,Du) & = \, \kabs{Du}^2 + {\rm Tr} (Du)^2 \, = \, 2 \, \kabs{x}^{-2} \,.
\end{align*}
From the first line we immediately obtain $u \in W^{1,p}(B,\R^2)$ and $\kabs{Du(x)} \geq 1$ for all $x \in B \setminus \{0\}$, which implies $g(\kabs{Du}^2) = 0$. In order to verify that $u$ a minimizer we start by recalling the identity $\diverg ( \kabs{x}^{-n-1} x \otimes x) = 0$ for all $x \in \R^n \setminus \{0\}$. Applying this in the two-dimensional case, we hence arrive at
\begin{align*}
\lefteqn{ \hspace{-0.5cm} \diverg \big[ \big( g(\kabs{Du}^2) + m_g \, A(x/\kabs{x})(Du,Du) \big)^{\frac{p-2}{2}} \, \big( g'(\kabs{Du}^2) \, Du + m_g \, A(x/\kabs{x}) Du \big) \big] } \\
	& = \, 2^{\frac{p-2}{2}} \, m_g^{\frac{p}{2}}\sum_{\kappa \in \{1,2\}} \frac{d}{dx_{\kappa}} 
	\Big[ \kabs{x}^{2-p} \, \big(A(x/\kabs{x}) Du\big)^{\kappa} \Big] \\
	&  = \, 2^{\frac{p-2}{2}} \, m_g^{\frac{p}{2}} \, \Big[ (2-p) \,  
	 \frac{x}{\kabs{x}^p} \, \Big( \frac{2}{\kabs{x}} + \frac{2(p-1)}{2-p} 
	\, \frac{x_1^2 + x_2^2}{\kabs{x}^3}\Big) 
	+  \kabs{x}^{2-p} \, 2 \, \diverg \frac{\II}{\kabs{x}} \Big] \\
	& = \, 2^{\frac{p-2}{2}} \, m_g^{\frac{p}{2}}  \, \Big[  
	(2-p) \, \frac{x}{\kabs{x}^{p+1}}
	\, \frac{2}{2-p} + \kabs{x}^{2-p}
	\, \frac{-2x}{\kabs{x}^3} \Big] \, = \, 0 \,.
\end{align*}
Therefore, observing that the expression in the divergence on the left-hand side of the previous equality belongs to $W^{1,1}(B,\R^{2 \times 2})$, we derive from the integration by parts formula that $u$ is a weak solution to the Euler-Lagrange system \eqref{EL-Ex} to \eqref{minimizer}, which means that $u$ is a critical point. The strict convexity of $f$ then yields the minimization property and the uniqueness, and this concludes the proof.
\end{proof}
\end{proposition}

In particular, the Euler-Lagrange system \eqref{EL-Ex} has a discontinuous solution (with coefficients still depending only on the independent and the gradient variable). Moreover, we may also take advantage of the particular structure of the integrand or the coefficients above, in the sense that the $x$-dependence occurs only in terms of $x/\kabs{x}$, i.\,e. of the minimizer itself. Expressing the $x$-dependence through the known solution (and omitting the $g'$-term which anyway vanishes for $u$) leads to the following definition of coefficients $a(u,z) \colon \R^2 \times \R^{2 \times 2} \to \R^{2 \times 2}$ by
\begin{equation}
\label{coeff}
a(u,z) \, := \, \big( g(\kabs{z}^2) + m_g \, A(u) (z,z) \big)^{\frac{p-2}{2}} \, A(u) z
\end{equation}
for all $(u,z) \in \R^2 \times \R^{2 \times 2}$. We deduce essentially from Proposition~\ref{prop-GV-f} that these coefficients have the correct behavior concerning growth and ellipticity, and we further prove a continuity assumption with respect to the $u$-variable:

\begin{proposition}
\label{prop-ass-coefficients}
The coefficients $a(\cdot,\cdot)$ defined in \eqref{coeff} are smooth with respect to the variable $z$ and satisfy the assumptions \eqref{GV} with constants $\nu,L$ depending only on $p,m_g$ and with $\alpha = 1$.

\begin{proof}
We first observe that the coefficients are smooth in the gradient variable by definition and choice of the cut-off function $g$. We now use $p<2$ and the boundedness of the bilinear form $A$ (by a constant depending only on $p$, independently of $u$) to find
\begin{equation*}
\kabs{a(u,z)} \, \leq \,  c(p,m_g) \, \big( 1 + \kabs{z} \big)^{p-1}\,,
\end{equation*}
where the constant $c(p,m_g)$ blows up as $p \nearrow 2$. Furthermore, the boundedness and ellipticity of $D_z a(u,z)$ is proved as in Proposition~\ref{prop-GV-f} (with the slight simplification that the second derivative of $g$ does not appear any more). Thus, the assumptions \eqref{GV}$_2$ and \eqref{GV}$_3$ hold true, and it only remains to verify the continuity condition \eqref{GV}$_4$: here we note that the bilinear form $A(u)$ is differentiable with respect to $u$ with bounded derivatives. Therefore, also $a(u,z)$ is differentiable with respect to $u$, and $a(u,z)$ and $D_u a(u,z)$ are bounded by $c(p,m_g) (1 + \kabs{z})^{p-1}$. Distinguishing the cases $\kabs{u - \bar{u}} \geq 1$ and $\kabs{u - \bar{u}} < 1$, we hence end up with
\begin{equation*}
\kabs{a(u,z) - a(\bar{u},z)} \, \leq \, c(p,m_g) \, \min\{\kabs{u - \bar{u}},1\} \, \big( 1 +  \kabs{z} \big)^{p-1}
\end{equation*}
for all $u,\bar{u} \in \R^2$ and all $z \in \R^{2 \times 2}$. Thus, the assumptions in \eqref{GV} are satisfied with the asserted dependencies.
\end{proof}
\end{proposition}

\begin{remark}
\label{rem-blow-up}
We emphasize that the ellipticity ratio $L/\nu$ of the coefficients $a(\cdot,\cdot)$ blows up as $p \nearrow 2$ by definition of the bilinear form $A$, and this property is indeed necessary for the construction of an elliptic system with a discontinuous weak solution in view of the existence of the ``critical'' exponent $p_1$ in Theorem~\ref{reg_systems} (the higher integrability exponent $q>p$ for $Du$  depends only on the structure data; in particular, the difference $q-p$ approaches zero as the ellipticity ratio $L/\nu$ blows up).
\end{remark}

It is now an easy consequence of Proposition~\ref{prop-min-Ex} that there exists a discontinuous weak solution -- namely the function $x/\kabs{x}$ as above -- to the homogeneous system related to the coefficients given by \eqref{coeff}. Hence, we have an example of a system satisfying the assumptions \eqref{GV} and admitting a weak solution with non-empty singular sets:

\begin{theorem}
\label{theorem-ex}
Assume $p \in (1,2)$ and let $u \colon \R^2 \supset B \to \R^2$ be given by $u(x) = x/\kabs{x}$. Then $u \in W^{1,p}(B,\R^2) \cap L^{\infty}(B,\R^2)$, and $u$ is a weak solution of the system
\begin{equation}
\label{system-example}
\diverg a(u,Du) \, = \, 0 \qquad \text{in } B \,,
\end{equation}
with coefficients $a(\cdot,\cdot)$ defined in \eqref{coeff}.

\begin{proof}
We observe that the choice $u(x)=x/\kabs{x}$ implies $a(u,Du) = D_z f(\cdot,Du)$ with $f$ taken from \eqref{integrand}. Thus the assertion follows from Proposition~\ref{prop-min-Ex}, where $\diverg D_z f(\cdot,Du) = 0$ was calculated.
\end{proof}
\end{theorem}

\begin{remarkon}
The theorem also provides an example that the fractional differentiability $Du \in W^{\alpha,p}$ cannot be obtained in the general subquadratic case for weak solutions to elliptic systems depending also explicitly on the solution itself.
\end{remarkon}

\begin{remarkon}
The question remains open whether or not there exists a discontinuous minimizer of a variational integral ${\cal F}[\cdot]$ with an integrand satisfying all the assumptions \eqref{GV-f}.  Instead of replacing $x/\kabs{x}$ in the coefficients in the Euler-Lagrange equation \eqref{EL-Ex} by $u$, we could also have argued on the level of the integrand \eqref{integrand}, defining
\begin{equation*}
\tilde{f}(u,z) \, := \, \big( g(\kabs{z}^2) +  m_g \, A(u)(z,z) \big)^{\frac{p}{2}}
\end{equation*}
for all $(u,z) \in \R^2 \times \R^{2 \times 2}$ and then studying minimizers of the associated variational integral. It is then easy to check that the function $x/\kabs{x}$ is still a critical point, but due to the lack of convexity of the integrand with respect to the $u$-variable this does not necessarily imply the minimization property  -- and hence it does not lead in a straightforward way to a counterexample to full regularity for minimizers. So far it is not clear if minimality holds (as in the quadratic case for $n\geq 3$ for which Giusti and Miranda were able to take advantage of the Euler-Lagrange equation) or not.
\end{remarkon}


\footnotesize

\begin{thebibliography}{Cam87b}

\bibitem[Ada75]{ADAMS75}
\textsc{Adams, R.A.}: \emph{{Sobolev Spaces}}.
\newblock Academic Press, New York (1975).

\bibitem[AF87]{ACEFUS87}
\textsc{Acerbi, E., Fusco, N.}: {A regularity theorem for minimizers of
  quasiconvex integrals}.
\newblock \emph{Arch. Ration. Mech. Anal.} {\bf 99} (1987), 261--281.

\bibitem[Bec09]{BECK09b}
\textsc{Beck, L.}: {Partial H\"older continuity for solutions of subquadratic
  elliptic systems in low dimensions}.
\newblock \emph{J. Math. Anal. Appl.} {\bf354} (2009), 1, 301--318.

\bibitem[Bec10]{BECK09}
\textsc{Beck, L.}: Boundary regularity results for variational integrals.
\newblock \emph{accepted for publication in Q. J. Math.}  (2010).

\bibitem[Cam82]{CAMPANATO82}
\textsc{Campanato, S.}: {H\"older continuity and partial H\"older continuity
  results for $W^{1,q}$-solutions of non-linear elliptic systems with
  controlled growth}.
\newblock \emph{Rend. Sem. Mat. Fis. Milano} {\bf 52} (1982), 435--472.

\bibitem[Cam87a]{CAMPANATO87b}
\textsc{Campanato, S.}: {A maximum principle for nonlinear elliptic systems:
  Boundary fundamental estimates}.
\newblock \emph{Adv. Math.} {\bf 66} (1987), 291--317.

\bibitem[Cam87b]{CAMPANATO87}
\textsc{Campanato, S.}: {Elliptic systems with non-linearity $q$ greater or
  equal to two. Regularity of the solution of the Dirichlet problem}.
\newblock \emph{Ann. Mat. Pura Appl. Ser. 4} {\bf 147} (1987), 117--150.

\bibitem[DG68]{DEGIORGI68}
\textsc{De~Giorgi, E.}: {Un esempio di estremali discontinue per un problema
  variazionale di tipo ellittico}.
\newblock \emph{Boll. Unione Mat. Ital., IV.} {\bf 1} (1968), 135--137.

\bibitem[DKM07]{DUZKRIMIN05}
\textsc{Duzaar, F., Kristensen, J., Mingione, G.}: {The existence of regular
  boundary points for non-linear elliptic systems}.
\newblock \emph{{J. Reine Angew. Math.}} {\bf 602} (2007), 17--58.

\bibitem[Eva86]{EVANS86}
\textsc{Evans, L.C.}: Quasiconvexity and partial regularity in the calculus of
  variations.
\newblock \emph{Arch. Ration. Mech. Anal.} {\bf95} (1986), 227--252.

\bibitem[FH85]{FUSHUT86}
\textsc{Fusco, N., Hutchinson, J.E.}: {$C^{1,\alpha}$ partial regularity of
  functions minimising quasiconvex integrals}.
\newblock \emph{Manuscr. Math.} {\bf 54} (1985), 121--143.

\bibitem[Gia83]{GIAQUINTA83}
\textsc{Giaquinta, M.}: \emph{{Multiple Integrals in the Calculus of Variations
  and Nonlinear Elliptic Systems}}.
\newblock Princeton University Press, Princeton, New Jersey (1983).

\bibitem[GM68a]{GIUSMIRA68}
\textsc{Giusti, E., Miranda, M.}: {Sulla Regolarit\`{a} delle Soluzioni Deboli
  di una Classe di Sistemi Ellitici Quasi-lineari}.
\newblock \emph{Arch. Rational Mech. Anal.} {\bf 31} (1968), 173--184.

\bibitem[GM68b]{GIUSMIRA68EX}
\textsc{Giusti, E., Miranda, M.}: {Un esempio di soluzioni discontinue per un
  problema di minimo relativo ad un integrale regolare del calcolo delle
  variazioni}.
\newblock \emph{Boll. Unione Mat. Ital., IV. Ser.} {\bf 1} (1968), 219--226.

\bibitem[GM79]{GIAMOD79}
\textsc{Giaquinta, M., Modica, G.}: {Almost-everywhere regularity results for
  solutions of non linear elliptic systems}.
\newblock \emph{Manuscr. Math.} {\bf 28} (1979), 109--158.

\bibitem[GM86]{GIAMOD86}
\textsc{Giaquinta, M., Modica, G.}: {Partial regularity of minimizers of
  quasiconvex integrals}.
\newblock \emph{{Ann. Inst. H. Poincar\'e Anal. Non Lin\'eaire}} {\bf 3}
  (1986), 185--208.

\bibitem[Gro02]{GROTOWSKI02}
\textsc{Grotowski, J.F.}: {Boundary regularity results for nonlinear elliptic
  systems}.
\newblock \emph{Calc. Var. Partial Differ. Equ.} {\bf 15} (2002), 353--388.

\bibitem[HLN96]{HAOLEONEC96}
\textsc{Hao, W., Leonardi, S., Ne\v{c}as, J.}: {An example of irregular
  solution to a nonlinear Euler-Lagrange elliptic system with real analytic
  coefficients}.
\newblock \emph{Ann. Sc. Norm. Super. Pisa, Cl. Sci., IV. Ser.} {\bf23} (1996),
  57--67.

\bibitem[HLW98]{HARLINWAN98}
\textsc{Hardt, R., Lin, F., Wang, C.}: {The $p$-energy minimality of $x/| x|$}.
\newblock \emph{Commun. Anal. Geom.} {\bf6} (1998), 1, 141--152.

\bibitem[Hon01]{HONG01}
\textsc{Hong, M.C.}: {On the minimality of the $p$-harmonic map
  $\frac{x}{|x|}:B^n\to S^{n-1}$}.
\newblock \emph{Calc. Var. Partial Differ. Equ.} {\bf13} (2001), 4, 459--468.

\bibitem[KM06]{KRIMIN05}
\textsc{Kristensen, J., Mingione, G.}: {The Singular Set of Minima of Integral
  Functionals}.
\newblock \emph{Arch. Rational Mech. Anal.} {\bf 180} (2006), 3, 331--398.

\bibitem[Min03a]{MINGIONE03}
\textsc{Mingione, G.}: Bounds for the singular set of solutions to non linear
  elliptic systems.
\newblock \emph{Calc. Var. Partial Differ. Equ.} {\bf 18} (2003), 4, 373--400.

\bibitem[Min03b]{MINGIONE03arch}
\textsc{Mingione, G.}: {The Singular Set of Solutions to Non-Differentiable
  Elliptic Systems}.
\newblock \emph{Arch. Rational Mech. Anal.} {\bf 166} (2003), 287--301.

\bibitem[Min06]{Mingione06Dark}
\textsc{Mingione, G.}: {Regularity of minima: an invitation to the Dark Side of
  the Calculus of Variations}.
\newblock \emph{Appl. Math.} {\bf 51} (2006), 4, 355--425.

\bibitem[Mor68]{MORREY68}
\textsc{Morrey, C.B.}: {Partial regularity results for non-linear elliptic
  systems}.
\newblock \emph{J. Math. Mech.} {\bf17} (1968), 649--670.

\bibitem[Nec77]{NECAS77}
\textsc{Necas, J.}: {Example of an irregular solution to a nonlinear elliptic
  system with analytic coefficients and conditions for regularity}.
\newblock {Theor. Nonlin. Oper., Constr. Aspects, Proc. int. Summer Sch.,
  Berlin 1975, 197--206} (1977).

\bibitem[Phi83]{PHILLIPS83}
\textsc{Phillips, D.}: A minimization problem and the regularity of solutions
  in the presence of a free boundary.
\newblock \emph{Indiana Univ. Math. J.} {\bf 32} (1983), 1--17.

\bibitem[SS10]{SCHEVENSCHMIDT10}
\textsc{Scheven, C., Schmidt, T.}: {Asymptotically regular problems I: Higher
  integrability}.
\newblock \emph{J. Differ. Equations} {\bf 248} (2010), 4, 745--791.

\bibitem[Sta71]{STARA71}
\textsc{Star\'a, J.}: Regularity results for non-linear elliptic systems in two
  dimensions.
\newblock \emph{Ann. Sc. Norm. Super. Pisa, Sci. Fis. Mat., III. Ser.} {\bf25}
  (1971), 163--190.

\bibitem[Wid71]{WIDMAN71}
\textsc{Widman, K.O.}: H\"older continuity of solutions of elliptic systems.
\newblock \emph{Manuscr. Math.} {\bf5} (1971), 299--308.

\end{thebibliography}
\end{document}